\documentclass[a4paper,12pt]{amsart}
\usepackage{ifthen}
\usepackage{graphicx}
\usepackage{mathrsfs}
\usepackage{color}
\nonstopmode \numberwithin{equation}{section}
\newtheorem{thm}{Theorem}[section]

\newtheorem{lem}[thm]{Lemma}

\newtheorem{prob}[thm]{Problem}

\theoremstyle{definition}

\newenvironment{pf}[1][]{%
 \vskip 3mm
 \noindent
 \ifthenelse{\equal{#1}{}}%
  {{\slshape Proof. }}%
  {{\slshape #1.} }%
 }%
{\qed\bigskip}

\newcounter{alphabet}
\newcounter{tmp}
\newenvironment{Thm}[1][]{\refstepcounter{alphabet}%
\bigskip%
\noindent%
{\bf Theorem \Alph{alphabet}}%
\ifthenelse{\equal{#1}{}}{}{ (#1)}%
{\bf .} \itshape}{\vskip 8pt}

\newcommand{\C}{{\mathbb C}}
\newcommand{\D}{{\mathbb D}}

\newcommand{\R}{{\mathbb R}}

\renewcommand{\Im}{{\,\operatorname{Im}\,}}

\renewcommand{\Re}{{\,\operatorname{Re}\,}}

\newcommand{\B}{{\operatorname{B}}}

\newcommand{\Gauss}{{\null_2F_1}}

\renewcommand{\arg}{\,{\operatorname{arg}\,}}

\newcommand{\aand}{{\quad\text{and}\quad}}

\newcounter{minutes}
\divide\time by 60
\newcounter{hours}
\multiply\time by 60 \addtocounter{minutes}{-\time}

\begin{document}
\bibliographystyle{amsplain}
\title{Mapping properties of the
zero-balanced hypergeometric functions}

\begin{center}
{\tiny \texttt{FILE:~\jobname .tex,
        printed: \number\year-\number\month-\number\day,
        \thehours.\ifnum\theminutes<10{0}\fi\theminutes}
}
\end{center}
\author[L.-M.~Wang]{Li-Mei Wang}
\address{School of Statistics,
University of International Business and Economics, No.~10, Huixin
Dongjie, Chaoyang District, Beijing 100029, China}
\email{wangmabel@163.com}


\begin{abstract}
In the present paper,
the order of convexity of $z\Gauss(a,b;c;z)$
is first given under some conditions on the
positive real parameters $a, b$ and $c$.
Then we show that
the image domains
of the unit disc $\D$ under some shifted zero-balanced
hypergeometric functions $z\Gauss(a,b;a+b;z)$
are convex and bounded by two horizontal lines which
solves the problem raised by Ponnusamy and Vuorinen
in \cite{PonVuor:2001}.

\end{abstract}

\thanks{This research is supported by 
National Natural Science Foundation of China (No.
11901086) " and the Fundamental
Research Funds for the Central Universities" in UIBE (No. 18YB02).
}
\maketitle

\section{Introduction and main results}

The
{\it Gaussian hypergeometric function}
 is defined by the power series
$$
\Gauss(a,b;c;z)=\sum_{n=0}^{\infty}\frac{(a)_n(b)_n}{(c)_nn!}z^n
$$
for $z\in\D:=\{z\in\C: |z|<1\},$
where $(a)_n$ is the {\it Pochhammer symbol}; namely, $(a)_0=1$ and
$(a)_{n+1}=(a)_n(a+n)=a(a+1)\cdots(a+n)$
for all $n\in\mathbb{N}:=\{0,1,2,\dots\}.$
Here $a,~b$ and $c$ are complex constants with
$-c\not\in\mathbb{N}$.
Hypergeometric functions can be
analytically continued along any path in the complex plane that avoids the
branch points $1$ and $\infty$.
If $c=a+b$, it is called the {\it zero-balanced} one
and $z\Gauss(a,b;c;z)$ is usually called a
{\it shifted} hypergeometric function.
For instance the function $z\Gauss(1, 1; 2; z)=-\log(1-z)$ 
is a shifted
zero-balanced hypergeometric function.
In the present paper, we only restrict to the real parameters
$a,~b$ and $c$.
For the basic properties of hypergeometric functions we refer to
\cite{AbramowitzStegun:1965}, \cite{OLBC:2010} and \cite{Wall:anal}.

For a function $f$ analytic in $\D$ and
normalized by $f(0)=f'(0)-1=0$,
the {\it order of convexity} and
the {\it order of starlikeness} (with
respect to the origin)
 of $f$ are defined by
 $$
\kappa=\kappa(f):=1+\inf_{z\in \D}\Re \frac{zf''(z)}{f'(z)}\in[-\infty,1]
$$
and
$$
\sigma=\sigma(f):=\inf_{z\in \D}\Re \frac{zf'(z)}{f(z)}\in[-\infty,1]
$$
respectively.
We follow the convention in \cite{Kustner:2007} due to K\"ustner,
that is
$\kappa(f)=-\infty$ only
if $f'$ is zero-free in $\D$ and $\Re [zf''(z)/f'(z)]$
is not bounded from below in $\D$, whereas
$\kappa(f)$ is regarded to be not defined if $f'$
has zeros in $\D$.
The corresponding convention can also be made for
$\sigma(f)=-\infty$.
It is known that $f$ is {\it convex}, i.e. $\kappa(f)\geq 0$
if and only if $f$ is univalent in $\D$ and $f(\D)$ is a convex domain;
and $f$ is {\it starlike}, i.e. $\sigma(f)\geq 0$
if and only if $f$ is univalent in $\D$ and $f(\D)$ is a starlike domain
with respect to the origin.
It is also true that if $\kappa(f)\geq -1/2$,
then $f$ is univalent in $\D$ and $f(\D)$ is convex in
(at least) one direction, see \cite{Umezawa:1952}
and \cite[p.17, Thm.2.24; p.73]{Rus:conv}.
If $f$ is convex in $\D$, then the order of
starlikeness of $f$ is at least $1/2$
(see \cite[p. 49, Theorem 2.1]{Rus:conv}).

The study on the order of starlikeness of the shifted
hypergeometric functions is dated from 1961,
see \cite{MS61} by Merkes and Scott.
Since then, several authors have researched
this project by using different methods,
for instance Ruscheweyh and Singh \cite{RuscheweyhSingh:1986},
Ponnusamy and Vuorinen \cite{PonVuor:2001},
K\"ustner \cite{Kustner:2002} and so on.
With regard to the order of convexity,
there are few papers in the literature
comparing with the order of starlikeness.
Such kind of research can be found in
\cite[Theorem 4]{Silv:1993} due to Silverman,
\cite[Corollary 8, Corollary 9]{Kustner:2007} by K\"ustner
in 2007 and our paper
\cite[Theorem 1.3, Theorem 1.5]{Wang:2020}.
Our first aim in this paper is to provide explicit
orders of convexity in terms of $a,~b$ and $c$
for some shifted hypergeometric functions.
As an application, we solve the following problem
on the zero-balanced hypergeometric functions
which is raised by Ponnusamy and Vuorinen in \cite{PonVuor:2001},
see also \cite{HPV:2010}.

\begin{prob}
\label{prob}
Do there exist positive numbers $\delta_1$, $\delta_2$ such that
 for $a\in (0,\delta_1)$ and $b\in(0,\delta_2)$
the normalized function $z\Gauss(a, b; a + b; z)$
$($$z\Gauss(a, b; a + b; z^2)$ resp.$)$
 satisfies the property that
maps the unit disc $\D$ into a strip domain?
\end{prob}

We restrict our attention
 only on the mapping $z\mapsto z\Gauss(a, b; a + b; z)$
in the present paper.
In view of Ramanujan's formula:
\begin{equation}\label{eq:Ramanujan}
\Gauss(a,b;a+b;z)=\frac{1}{\B(a,b)}\big(R(a,b)-\log(1-z)\big)
+O\left(|1-z|\log\frac1{|1-z|}\right),
\end{equation}
as $z\to1$ in $\D,$ where
$
R(a,b)=2\psi(1)-\psi(a)-\psi(b)
$
and $\psi(x)=\Gamma'(x)/\Gamma(x)$
and $\B(a,b)$ denote the digamma function and the Beta
function respectively,
the zero-balanced hypergeometric function behaves
like the function $-\frac{1}{\B(a,b)}\log(1-z)$. 
Thus it is easy to see that
it maps the unit disc $\D$ into a strip domain.
It seems that it is trivial to consider only
Problem 1.1 along this direction.
In fact, our aim is to find conditions on the
pair $(a,b)$ such that the image domain of
$\D$ under the mapping $z\mapsto z\Gauss(a, b; a + b; z)$
is convex and lies in a parallel strip.
Along this line,
the author in \cite[Corollary 1.2]{Wang:2020}
(cf. \cite[Corollary 8, Corollary 9]{Kustner:2007})
got the following result:

\begin{Thm}
\label{Wang1}
If $a$ and $b$ are real constants
satisfying $-a\not\in\mathbb{N}$, $-b\not\in\mathbb{N}$
and $ab<1$, then
the shifted zero-balanced hypergeometric function
$z\Gauss(a, b; a+b; z)$ is not convex in $\D$.
\end{Thm}

If furthermore, we suppose $a$ and $b$ are 
both positive. According to Theorem A,
in order to consider the
convexity of the shifted zero-balanced
hypergeometric functions the condition
$ab\geq 1$ is required,
which further implies $a\geq 1$
or $b\geq 1$.
Furthermore if $a+b=2$,
it follows from $2=a+b\geq 2\sqrt{ab}$
that $ab\leq 1$
which together with $ab\geq 1$
implies $ab=1$.
Thus when investigating the
convexity of the shifted zero-balanced
hypergeometric functions
$z\Gauss(a,2-a;2;z)$ with $a>0$, it is sufficient
to consider only the function
$z\Gauss(1,1;2;z)=-\log(1-z)$.

For $0<a\leq b<1$, even though $z\Gauss(a, b; a+b; z)$
is not convex,
it is natural to ask the question:
whether the function
$z\Gauss(a, b; a+b; z)$ possesses the other
geometric properties,
for example starlikeness and univalence
and so on.
In fact,  K\"ustner in \cite[Theorem 1]{Kustner:2007}
(cf. \cite[Theorem 1.1, Remark 2.3]{Kustner:2002})
got the following result on the order of starlikenss of
the shifted hypergeometric functions.

\begin{Thm}\label{thm:starlike}
If $0<a\leq b\leq c$, then
$$
\sigma(z\Gauss(a, b; c; z))=1-\frac{\Gauss'(a, b; c; -1)}{\Gauss(a, b; c; -1)}
\geq 1-\frac{ab}{b+c}.
$$
\end{Thm}
As a consequence of Theorem B,
the function $z\Gauss(a, b; a+b; z)$
is starlike thus univalent in $\D$
if $0<a\leq b<1$.

From now on, we focus our attention on
the case
$a\geq 1$. First if $a=1$,
we proved the next result in
\cite{Wang:2020} on the convexity of
$z\Gauss(1, b; c; z)$.

\begin{Thm}$($\cite[Corollary 4.3]{Wang:2020}$)$
\label{Wang2}
Assume $b$ and $c$ are real
constants satisfying $0<b\leq c$,
the function $z\Gauss(1, b; c; z)$
is convex in $\D$
if one of the following
conditions holds:
\begin{enumerate}
\item $0\leq b\leq 4$ and $c\geq 2$;
\item $1\leq c<\min\{2,1+b\}$
and $c\geq 3-b$.
\end{enumerate}
\end{Thm}

For general $a>1$, by making use of Riemman-Stieltjes type
integral of the ratio of hypergeometric functions,
we first obtain the order of convexity of
the shifted hypergeometric functions
as follows.

\begin{thm}\label{thm:order}
Let $a, b$
and $c$ be real constants
satisfying $1<a\leq b\leq 4$, $c>2$ and
$b\leq c\leq 2b\leq 2(1+a)$.
The order of convexity of the function $z\Gauss(a, b; c; z)$
is
\begin{eqnarray*}
\kappa(z\Gauss(a, b; c; z))
&=&\frac{5-c-a-b}{2}+\frac{c-2-(1-a)(1-b)}{2\left(1-a+a
\dfrac{\Gauss(a+1, b; c; -1)}{\Gauss(a, b; c; -1)}\right)}.\\
\end{eqnarray*}
In particular, $\Gauss'(a, b; c; z)\not=0$
for all $z\in\D$.
\end{thm}

We remark that in \cite{Wang:2020},
the order of convexity of $z\Gauss(a, b; c; z)$
has already been dealt with, but the condition
$0<a\leq 1$ is required there. Thus Theorem
\ref{thm:order} extends
the discussion to $a>1$.
By virtue of Theorems C and \ref{thm:order},
we obtain the next theorem which not only
gives the convexity of some shifted
zero-balanced hypergeometric functions
but also
demonstrates the explicit bounds of
the strip domains of the image domains.
Thus we provide a strong version
of solution to Problem \ref{prob}.

\begin{thm}\label{thm:conv}
Let $a$ and $b$ be real constants.
If one of the following conditions holds:
\begin{enumerate}
\item $a=1$ and $1\leq b\leq 4$;
\item $1<a\leq b\leq \min\{3,1+a,\frac{4a}{5a-4}\}$,
\end{enumerate}
then
the shifted zero-balanced hypergeometric function
$z\Gauss(a, b; a + b; z)$
is convex and
maps the unit disc $\D$ into a strip domain
which is bounded by two horizontal lines
$\Im w=\frac{\pi}{2\B(a,b)}$.
The boundary lines are optimal. 
\end{thm}


\section{Some lemmas}
In this section, we list some auxiliary lemmas
which play important roles in the proofs of
the main results.
The first lemma is due to Ruscheweyh,
Salinas and Sugawa \cite{RusSalSu}.
Later, Liu and Pego in \cite{LiuPego:2016}
pointed out that one condition can be deduced
from the others.
So the lemma can be finally
characterized in the following form.

\begin{lem}\label{measure}
Let $F(z)$ be analytic in the slit domain $\C\setminus[1,+\infty)$. Then
$$
F(z)=\int_{0}^{1}\frac{d\mu(t)}{1-tz}
$$
for some probability measure $\mu$ on $[0,1]$, if and only if the following conditions are fulfilled:
\begin{enumerate}
\item $F(0)=1$;\\
\item $F(x)\in\R$ for $x\in(-\infty,1)$;\\
\item $\Im F(z)\ge 0$ for $\Im z>0$;\\
\item $\displaystyle \limsup_{x\to+\infty}F(-x)\ge 0$.
\end{enumerate}
The measure $\mu$ and the functions $F$ are in one-to-one correspondence.
\end{lem}

The forthcoming four lemmas describe the
properties of the ratio to two hypergeometric functions
in different aspects.
\begin{lem}
$($\cite[Thm. 1.5]{Kustner:2002},  \cite[ p.337-339 and Thm.69.2]{Wall:anal}$)$
\label{lem:integral}
If $-1\leq a\leq c$ and $0\leq b\leq c\not=0$, the ratio of two
hypergeometric functions can be written in
integral as
$$
\frac{\Gauss(a+1,b+1;c+1;z)}{\Gauss(a+1,b;c;z)}
=\int_{0}^{1}\frac{d\mu(t)}{1-tz},
\quad z\in\C\setminus[1,+\infty)
$$
where $\mu:[0,1]\to[0,1]$ is nondecreasing with $\mu(1)-\mu(0)=1$.
\end{lem}

\begin{lem}\label{lem:z=-infty}
Let $a,~b$ and $c$ be real constants with
$0<a\leq b\leq 1+a$, $c-a\not\in-\mathbb{N}$
and $c-b\not\in-\mathbb{N}$.
 Then
$$
\displaystyle\lim_{x\to+\infty}
\frac{x\Gauss(a+1,b+1;c+1;-x)}{\Gauss(a+1,b;c;-x)}=\infty.
$$
\end{lem}
\begin{pf}
Denote $H(z)=\Gauss(a+1,b+1;c+1;z)$
and $G(z)=\Gauss(a+1,b;c;z)$
for convenience.
In order to show the claimed equation,
the next two linear transforms
(see \cite[p.559, 15.3.3, 15.3.5, 15.3.7]{AbramowitzStegun:1965})
are required:
\begin{eqnarray}\label{eq:T1}
\notag\Gauss(a,b;c;z)
&=&\frac{\Gamma(c)\Gamma(b-a)}{\Gamma(b)\Gamma(c-a)}(-z)^{-a}\Gauss(a,1-c+a;1-b+a;1/z)\\
&+&\frac{\Gamma(c)\Gamma(a-b)}{\Gamma(a)\Gamma(c-b)}(-z)^{-b}\Gauss(b,1-c+b;1-a+b;1/z),
 \quad (|\arg(-z)|<\pi)
\end{eqnarray}
and
\begin{equation}\label{eq:T2}
\Gauss(a,b;c;z)=(1-z)^{-b}\Gauss\left(b,c-a;c;\frac{z}{z-1}\right),\quad
(|\arg(1-z)|<\pi).
\end{equation}
Then the proof can be separated
into three cases according to
the relationship between
$a$ and $b$.

Case I: Assume that $a<b<1+a$.
By virtue of \eqref{eq:T1},
we obtain
\begin{equation*}
H(-x)=
   \frac{\Gamma(c+1)\Gamma(b-a)}{\Gamma(b+1)\Gamma(c-a)}x^{-a-1}
   +O(x^{-b-1})
   \end{equation*}
and
\begin{equation*}
 G(-x)=
   \frac{\Gamma(c)\Gamma(a+1-b)}{\Gamma(a+1)\Gamma(c-b)}x^{-b}
   +O(x^{-a-1})
\end{equation*}
as $x\to+\infty$.
It follows from the above asymptotic behaviors
that
$
\displaystyle\lim_{x\to+\infty}
xH(-x)/G(-x)=\infty
$
since $a<b<1+a$.

Case II: Assume $b=a$. We deduce from the identity
\eqref{eq:T2} and Ramanujan's formula
\eqref{eq:Ramanujan} that
\begin{eqnarray*}
H(-x)&=&\Gauss(a+1,a+1;c+1;-x)
=(1+x)^{-a-1}\Gauss\left(a+1,c-a-1;c;\frac{x}{x+1}\right)\\
&=&(1+x)^{-a-1}\left(\frac{\log(1+x)}{\B(a+1,c-a-1)}+O(1)\right),
\end{eqnarray*}
as $x\to+\infty$.
Similarly, applying the identity
\eqref{eq:T2} and Gauss summation formula:
$$
\displaystyle\lim_{x\to 1^{-}}\Gauss(a,b;c;x)
=\Gauss(a,b;c;1)
=\frac{\Gamma(c)\Gamma(c-a-b)}{\Gamma(c-a)\Gamma(c-b)},
$$
if $c-a-b>0$ (see \cite[p.556, 15.1.20]{AbramowitzStegun:1965}),
 we have
\begin{eqnarray*}
\lim_{x\to+\infty}(1+x)^{a}G(-x)
&=& \lim_{x\to+\infty}\Gauss\left(a,c-a-1;c;\frac{x}{x+1}\right)
=\frac{\Gamma(c)}{\Gamma(c-a)\Gamma(a+1)}.
\end{eqnarray*}
Thus by substituting the proceeding descriptions of $H(-x)$
and $G(-x)$ as $x\to+\infty$, we deduce that
\begin{eqnarray*}
\displaystyle\lim_{x\to+\infty}
\frac{xH(-x)}{G(-x)}
&=&\lim_{x\to+\infty}\frac{x(1+x)^{-a-1}\left(\frac{\log(1+x)}{\B(a+1,c-a-1)}+O(1)\right)}{G(-x)}\\
&=&\lim_{x\to+\infty}\frac{x(1+x)^{-1}\left(\frac{\log(1+x)}{\B(a+1,c-a-1)}+O(1)\right)}{(1+x)^{a}G(-x)}
=\infty
\end{eqnarray*}
in this case.

Case III: Assume that $b=1+a$. Applying
the same techniques in Case II, we find that
$$
\lim_{x\to+\infty}(1+x)^{a+1}H(-x)=
\frac{\Gamma(c+1)}{\Gamma(c-a)\Gamma(a+2)},
$$
and
$$
G(-x)=(1+x)^{-a-1}\left(\frac{\log(1+x)}{\B(a+1,c-a-2)}+O(1)\right),
\quad x\to+\infty,
$$
which yields the required claim again.

We complete all these three cases,
thus the proof is done.
\end{pf}

\begin{lem}\label{lem:G/F}
Let $a, b$ and $c$ be real constants
with $a,b,c\not\in-\mathbb{N}$,
$c-a\not\in-\mathbb{N}$ and $c-b\not\in-\mathbb{N}.$
\begin{enumerate}
\item
If $c<a+b$, then
\begin{equation*}\label{beha3}
\frac{\Gauss(a+1,b;c;z)}{\Gauss(a,b;c;z)}
=\frac{a+b-c}{a(1-z)}+O\left(|1-z|^{a+b-c}\right).
\end{equation*}
\item
 If $c=a+b$, then
 \begin{equation*}\label{beha2}
\frac{\Gauss(a+1,b;c;z)}{\Gauss(a,b;c;z)}
=\frac{1}{-a(1-z)\log(1-z)}+O\left(\log\frac{1}{|1-z|}\right).
\end{equation*}
\item
If $a+b<c<a+b+1$, then
\begin{equation*}\label{beha1}
\frac{\Gauss(a+1,b;c;z)}{\Gauss(a,b;c;z)}
=\frac{A}{(1-z)^{1-\alpha}}+O(|1-z|^{\varepsilon-1})
\end{equation*}
where
$$
A=\frac{\Gamma(a+b+1-c)\Gamma(c-a)\Gamma(c-b)}
{\Gamma(a+1)\Gamma(b)\Gamma(c-a-b)},
$$
$\alpha=c-a-b\in(0,1)$ and $\varepsilon=\min\{2\alpha,1\}.$
\item
 If $c=1+a+b$, then
  \begin{equation}\label{beha4}
\frac{\Gauss(a+1,b;c;z)}{\Gauss(a,b;c;z)}
=-b\log(1-z)+O\left(1\right).
\end{equation}

\end{enumerate}
\end{lem}
\begin{pf}
The first three assertions can be found in
\cite[Lemma 2.3]{Wang:2020}. We need only to
prove the last one.
First recall the linear transform
\begin{align}
&\Gauss(a,b;c;z)
\label{tran}
=\frac{\Gamma(c)\Gamma(c-a-b)}{\Gamma(c-a)\Gamma(c-b)}
\Gauss(a,b;a+b-c+1;1-z) \\
&\qquad +(1-z)^{c-a-b}
\frac{\Gamma(c)\Gamma(a+b-c)}{\Gamma(a)\Gamma(b)}
\Gauss(c-a,c-b;c-a-b+1;1-z).
\notag
\end{align}
Applying the identity above and Ramanujan's
formula \eqref{eq:Ramanujan}, we have
\begin{eqnarray*}
&&\frac{\Gauss(a+1,b;c;z)}{\Gauss(a,b;c;z)}\\
&=&\frac{\frac{1}{\B(1+a,b)}\big(R(1+a,b)-\log(1-z)\big)
+O\left(|1-z|\log\frac1{|1-z|}\right)}{\frac{\Gamma(c)\Gamma(c-a-b)}
{\Gamma(c-a)\Gamma(c-b)}+O(|1-z|)}\\
&=&-b\log(1-z)+O(1)
\end{eqnarray*}
as $z\to1$ in $\D,$  since $c=1+a+b$.
\end{pf}

\begin{lem}$($\cite{Wang:2020},
see also \cite[Ramark 2.3]{Kustner:2002}$)$
\label{z=-1}
If $-1\leq a\leq c$ and $0\leq b\leq c\not=0$, then
$$
\frac{c}{b+c}\leq \frac{\Gauss(a+1, b; c; -1)}{\Gauss(a, b; c; -1)}
\leq \frac{2c-b}{2c}.
$$
\end{lem}

\begin{lem}$($see \cite[Lemma 3.2]{SugawaWang:2016}$)$
\label{lem:conv}
Let $\Omega$ be an unbounded convex domain in $\C$ whose
boundary is parametrized positively by a Jordan curve $w(t)=u(t)+iv(t),~
0<t<1,$ with $w(0^+)=w(1^-)=\infty.$
Suppose that $u(0^+)=+\infty$ and that $v(t)$ has a finite limit
as $t\to0^+.$
Then $v(t)\le v(0^+)$ for $0<t<1.$
\end{lem}

In the end of this section,
the behavior of a special zero-balanced
hypergeometric function $\Gauss(1,1,2;z)$
around $z=1$ is given.
\begin{lem}\label{lem:log}
For $\theta\in(0,2\pi)$,
we have
$$
\lim_{\theta\to0^+}\Re[-e^{i\theta}\log(1-e^{i\theta})]
=+\infty \aand
\lim_{\theta\to0^+}\Im[-e^{i\theta}\log(1-e^{i\theta})]
=\frac{\pi}{2}.
$$
\end{lem}

\begin{pf}
For $z=e^{i\theta}$ with $\theta\in(0,2\pi)$,
we do the following
computations:
\begin{eqnarray*}
&&-e^{i\theta}\log(1-e^{i\theta})
=-e^{i\theta}\left[2\log\left(\sin\frac{\theta}{2}\right)+i\frac{\theta-\pi}{2}\right]\\
&=&-2\cos\theta\log\left(\sin\frac{\theta}{2}\right)+\frac{\theta-\pi}{2}\sin\theta
-i\left[2\sin\theta\log\left(\sin\frac{\theta}{2}\right)+\frac{\theta-\pi}{2}\cos\theta\right].
\end{eqnarray*}
We deduce the required identities by letting
$\theta\to 0^+$ in the above form.
\end{pf}

%
%
%
%
%


\section{Proofs of the main results}

\begin{pf}[Proof of Theorem \ref{thm:order}]
Let $F(z)=\Gauss(a,b;c;z)$, $G(z)=\Gauss(a+1,b;c;z)$
and $H(z)=\Gauss(a+1,b+1;c+1;z)$
for simplicity.
In order to obtain the order of convexity of $zF(z)$,
we should first compute its second logarithmic
derivative and then evaluate the real part of
this derivative in the unit disc $\D$.
According to the equation (3.2) in \cite{Wang:2020},
we get
\begin{equation}\label{eq:W}
W(z):=1+\frac{z(zF)''}{(zF)'}=\frac{3-c+(a+b-2)z}{1-z}
+\frac{c-2+(1-a)(1-b)z}{(1-z)(1-a+aG/F)}.
\end{equation}
We will show that the infimum real part of $W(z)$
in $\D$ is attained at $z=-1$
under the assumptions of this theorem.

We infer from the linear transformation
(see \cite[eq. 27]{Gauss} or \cite[eq. (2.4)]{Kustner:2002})
$$
\dfrac{\Gauss(a+1,b;c;z)}{\Gauss(a,b;c;z)}
=\frac{1}{1-\dfrac{b}{c}\dfrac{z\Gauss(a+1,b+1;c+1;z)}
{\Gauss(a,b;c;z)}},
\quad z\in\C\setminus[1,+\infty)
$$
that
the denominator of the second term
in $W(z)$ can be transformed into
\begin{eqnarray*}
&&\frac{1}{1-a+aG/F}\\
&=&\frac{1}{1-a+a/[1-b/czH(z)/G(z)]}\\
&=&\frac{1}{1-a}+\frac{a}{a-1}
      \frac{1}{1+\frac{(a-1)b}{c}\frac{zH}{G}}.
\end{eqnarray*}
Substituting the equation above into equation (\ref{eq:W}),
we have
\begin{eqnarray*}
W(z)&=&\frac{3-c+(a+b-2)z}{1-z}+
\frac{c-2+(1-a)(1-b)z}{(1-a)(1-z)}\\
&+&\frac{a}{(a-1)(1-z)}
      \frac{c-2+(1-a)(1-b)z}{1+\frac{(a-1)b}{c}\frac{zH}{G}}\\
&=&1-a+\frac{a(a+1-c)}{a-1}\frac{1}{1-z}+
\frac{a}{(a-1)(1-z)}
      \frac{c-2+(1-a)(1-b)z}{1+\frac{(a-1)b}{c}\frac{zH}{G}}.
\end{eqnarray*}
Denote
$$
M_1(z)=\frac{c-2+(1-a)(1-b)z}{1-z},\quad
M_2(z)=1+\tau \frac{zH}{G}
$$
and $M(z)=\frac{M_1(z)}{M_2(z)}$ for simplicity
where $\tau=\frac{(a-1)b}{c}$.

Next we show that under the conditions of
this theorem, there exists a probability measure
$\nu$ on $[0,1]$ such that
\begin{equation}\label{eq:M}
M(z)
=(c-2)\int_{0}^{1}\frac{d\nu(t)}{1-tz}.
\end{equation}
By virtue of Lemma \ref{measure},
it is sufficient to verify the four
conditions there.

The first and second ones are easy to check
since $a,~b$ and $c$ are real constants with
$c>2$.
We proceed to verify the third one.
It follows from Lemma \ref{lem:integral}
that there exists a probability measures $\mu$
on $[0,1]$ such that
$$
\frac{H(z)}{G(z)}=\int_{0}^{1}\frac{d\mu(t)}{1-tz}
$$
which implies that
$$
M_2(z)=\int_{0}^{1}\frac{1+(\tau-t)z}{1-tz}d\mu(t).
$$
Thus for $z=x+iy$ with $y>0$,
an elementary computation generates that
\begin{eqnarray*}
&&\Im M(z)|M_2(z)|^2\\
&=&-\Re M_1\Im M_2+\Im M_1\Re M_2\\
&=&-\int_{0}^{1}\frac{\tau y}{|1-tz|^2}
\frac{c-2-(1-a)(1-b)|z|^2-[c-2-(1-a)(1-b)]x}{|1-z|^2}d\mu(t)\\
&&+\int_{0}^{1}
\frac{1+(t^2-\tau t)|z|^2+(\tau-2t)x}{|1-tz|^2}
\frac{[c-2+(1-a)(1-b)]y}{|1-z|^2}d\mu(t)\\
&:=&y\int_{0}^{1}\frac{I(t,x,y)}{|1-z|^2|1-tz|^2}d\mu(t).
\end{eqnarray*}

Letting $p=c-2+(1-a)(1-b)$
and $q=c-2-(1-a)(1-b)$ for simplicity,
the numerator of the integrand can be simplified
into
\begin{eqnarray*}
&&I(t,x,y)\\
&=&p(t^2-\tau t)|z|^2+\tau(1-a)(1-b)|z|^2+p(\tau -2t)x+\tau qx
+p-\tau(c-2)\\
&=&[\tau(1-a)(1-b)+p(t^2-\tau t)]|z|^2+[p(\tau-2t)+\tau q]x+p-\tau(c-2)\\
&:=&Q(t)y^2+Q(t)\left(x+\frac{\tau(c-2)-pt}{Q(t)}\right)^2+
  \frac{[p-\tau (c-2)]Q(t)-[\tau(c-2)-pt]^2}{Q(t)}\\
&:=&Q(t)y^2+Q(t)\left(x+\frac{\tau(c-2)-pt}{Q(t)}\right)^2+
  \frac{S(t)}{Q(t)},
\end{eqnarray*}
since $Q(t)\not\equiv0$
as $1<a\leq b$.
Thus for $y=\Im z>0$, $\Im M(z)\geq 0$ is valid if
\begin{equation*}
\begin{cases}
Q(t)=\tau(1-a)(1-b)+p(t^2-\tau t)\geq 0;\\
S(t)=[p-\tau (c-2)]Q(t)-[\tau(c-2)-pt]^2\geq 0,
\end{cases}
\end{equation*}
hold for any $t\in[0,1]$.

We prove the first inequality.
By an elementary calculation, we find that
\begin{equation}\label{inq:Q}
\begin{cases}
Q(0)=\tau(1-a)(1-b)>0;\\
Q(1)=c-2+\frac{(a-1)(2b-c)}{c}\geq 0.
\end{cases}
\end{equation}
Note also that the conditions $1<a\leq b\leq 4$
and $b\leq c$
imply $0<\tau\leq 3$.
Thus if $2\leq \tau\leq 3$, the equations in
\eqref{inq:Q} yield that the quadratic
function $Q(t)$ is nonnegative in
$[0,1]$.
Hence we remain to consider the case $0<\tau <2$,
and in this case the inequality
$$
Q(t)
\geq Q\left(\frac{\tau}{2}\right)
=\frac{\tau^2}{4b}[4(b-1)c-pb],
$$
always holds for $t\in[0,1]$.
Therefore in order to show the non-negativity
of $Q(t)$ on $[0,1]$, it suffices to prove
$4(b-1)c-pb\geq 0$.
As $1<a\leq b$ implies
$$
4(b-1)c-pb=(3b-4)c+2b-(a-1)(b-1)b\geq (3b-4)c+2b-(b-1)^2b
:=h(b,c),
$$
we divide the proof into two cases according to
the sign of $3b-4$.

Case I: Assume that $1<b\leq 4/3$. Then
$h(b,c)\geq 2b(3b-4)+2b-(b-1)^2b=-b(b-1)(b-7)
> 0$
as $c\leq 2b$.

Case II: Assume that $b>4/3$. Then
$h(b,c)\geq (3b-4)b+2b-(b-1)^2b=-b(b^2-5b+3)
\geq 0$
as $b\leq c$ and $4/3<b\leq 4$.

Thus we obtain $Q(t)\geq 0$ for $t\in[0,1]$ under
the assumptions of this theorem.

We next verify the property of $S(t)$.
After computing, we get
$$
S(t)=-p\tau[(c-2)t+p-(c-2)(\tau+2)]t+
[p-\tau(c-2)]\tau(1-a)(1-b)-(c-2)^2\tau^2
$$
and $S(1)=0$.
%
Furthermore an elementary computation yields that
$S(0)=pc\tau^2(2b-c)$, which is nonnegative as
$c\leq 2b$. Therefore
the quadratic function
$S(t)\geq0$ for all $t\in[0,1]$.

In the end, we turn to show the final condition
$\displaystyle \limsup_{x\to+\infty}M(-x)\ge 0$ in
Lemma \ref{measure}.
Lemma \ref{lem:z=-infty} shows that
$
\displaystyle\lim_{x\to+\infty}M_2(-x)=\infty
$
which in conjunction with
$
\displaystyle\lim_{x\to+\infty}M_1(-x)=-(a-1)(b-1)\not=0
$
generates that
$$
\lim_{x\to+\infty}M(-x)=0
$$
under the condition $a\leq b\leq 1+a$.

We thus obtain the integral form \eqref{eq:M} which immediately
yields that
$$
W(z)=1-a+\frac{a(a+1-c)}{a-1}\frac{1}{1-z}+
\frac{a(c-2)}{a-1}
\int_{0}^{1}\frac{ d\nu(t)}{1-tz},
$$
where $\nu$ is a probability measure on $[0,1]$.
Therefore $W(z)$
is analytic on $\D$ and
since $a>1$ and $c>2$, then
\begin{equation}\label{ineq:W}
\Re W(z)\geq W(-1),\quad \text{for}~|z|=1.
\end{equation}
Next we will prove the above
inequality for all $z\in \D$
which is exactly the assertion of this theorem.
If $c\leq a+1$ in addition,
it is obvious that inequality
\eqref{ineq:W} holds for $z\in\D$
since $\Re 1/(1-z)>1/2$ is valid
for $z\in\D$.
We need only to deal with the
case $c>a+1$.
First we observe that $c\leq 1+a+b$ as $c\leq 2b$
and $b\leq 1+a$.
As a direct consequence of Lemma \ref{lem:G/F},
we conclude that
\begin{equation*}\label{eq:M(1)}
\mathop{\lim_{z\to1^{-}}}_{z\in\R}
\frac{G(z)}{F(z)}=+\infty,
\end{equation*}
for $1+a<c\leq 1+a+b$.
Thus it follows from the above equation that
\begin{eqnarray*}
W(z)&=&
\frac{3-c+(a+b-2)z}{1-z}
+\frac{c-2+(1-a)(1-b)z}
{(1-z)(1-a+aG/F)}\\
&=&
\frac{3-c+(a+b-2)z}{1-z}
+\frac{o(1)}
{1-z}\\
&=&\frac{1+a+b-c+o(1)}{1-z}
\end{eqnarray*}
as $z\to 1$ along the real axis in $\D$.
Since $W(z)$ is analytic in $\D$,
we obtain that
\eqref{ineq:W} holds for all $z\in\D$
if $c<1+a+b$.
As for the case $c=1+a+b$, by virtue of equation
\eqref{beha4} in Lemma \ref{lem:G/F},
we conclude that
\begin{eqnarray*}
W(z)&=&
2-a-b
+\frac{c-2+(1-a)(1-b)z}
{1-z}\frac{1}{1-a-ab\log(1-z)+O(1)}\\
&\to&+\infty,
\end{eqnarray*}
as $z\to 1$ along the real axis in $\D$.
Thus
the inequality
\eqref{ineq:W} holds for all $z\in\D$ if $c=1+a+b$,
since $W(z)$ is analytic in $\D$.

We complete the proof.
\end{pf}

\begin{pf}[Proof of Theorem \ref{thm:conv}]
We first prove that $z\Gauss(a,b;a+b;z)$ is convex
in $\D$ under the assumptions. If $a=1$ and $1\leq b\leq 4$,
the convexity of $z\Gauss(a,b;a+b;z)$ is a
consequence of the case (2) in
Theorem C.
While if $1<a\leq b\leq 3$, the order of convexity of
$z\Gauss(a,b;a+b;z)$ is
\begin{equation*}
\kappa=\frac{5-2a-2b}{2}+\frac{2(a+b)-3-ab}{2\left(1-a+a
\dfrac{\Gauss(a+1, b; a+b; -1)}{\Gauss(a, b; a+b; -1)}\right)}
\end{equation*}
according to Theorem \ref{thm:order}.
It is easy to see that
$2(a+b)-3-ab\geq 0$ since $1<a\leq b\leq 3$.
Thus in view of Lemma \ref{z=-1},
we get
$$
\kappa\geq
\frac{5-2a-2b}{2}+\frac{2(a+b)-3-ab}{2\left(1-a+a
\frac{2a+b}{2(a+b)}\right)}
= \frac{(4-5a)b+4a}{2(2a+2b-ab)}.
$$
By observing the numerator of the above lower bound,
we conclude that if $1<a\leq b\leq 4a/(5a-4)$,
then the order of convexity $\kappa$ is non-negative,
which means $z\Gauss(a,b;a+b;z)$ is convex in $\D$.

On the other hand,
we consider the image domain $\Omega$ of $\D$
under the mapping $z\mapsto z\Gauss(a,b;a+b;z)$.
For $0<t<1$, let $w(t)=u(t)+iv(t)$ stand for
the boundary
$e^{2\pi i t }\Gauss(a,b;a+b;e^{2\pi i t })$
of the domain $G$, which is a Jordan convex curve.
In view of Lemma \ref{lem:log} and Ramanujan's formula
\eqref{eq:T2},
we have
$\displaystyle\lim_{t\to0^+}w(t)=
\displaystyle\lim_{t\to1^-}w(t)=\infty$,
$\displaystyle\lim_{t\to0^+}u(t)=+\infty$
and
$$
\lim_{t\to0^+}v(t)=
 \frac{\pi}{2\B(a,b)}.
$$
Therefore it follows from Lemma \ref{lem:conv}
that
$v(t)\leq v(0^+)$ holds for $0<t<1$.
Since the image domain $G$ is symmetric
with respect to the real axis, we finally
have
$$
| v(t)|\leq \frac{\pi}{2\B(a,b)},
\quad \text{for}~0<t<1,
$$
which means that the image domain of $\D$ under
the function $z\Gauss(a,b;a+b;z)$ is bounded by two
horizontal lines
$\Im w=\pm\frac{\pi}{2\B(a,b)}$.

The proof is finished.
\end{pf}

\noindent
{\bf Acknowledgements.}
The author would like to
appreciate Professor
Toshiyuki Sugawa for discussions and
guidance, without whom I could't finish
this work. The author is also indebted to
 Professor Matti
 Vuorinen for drawing
  my attention to this problem.

\def\cprime{$'$} \def\cprime{$'$} \def\cprime{$'$}
\providecommand{\bysame}{\leavevmode\hbox
to3em{\hrulefill}\thinspace}
\providecommand{\MR}{\relax\ifhmode\unskip\space\fi MR }
\providecommand{\MRhref}[2]{%
  \href{http://www.ams.org/mathscinet-getitem?mr=#1}{#2}
} \providecommand{\href}[2]{#2}


\begin{thebibliography}{10}

\bibitem{AbramowitzStegun:1965}
M. Abramowitz and I. Stegun, \emph{Handbook of Mathematical
Functions}, Dover, New York, 1965.





%
\bibitem{Gauss}
C. F. Gauss, Disquisitiones generales circa seriem infinitam... , Commentationes Societatis
Regiae Scientiarum Gottingensis Recentiores \textbf{2} (1812), 1-46; reprint in C. F. Gauss, Werke,
Band III, K\"{o}onigliche Gesellschaft der Wissenschaften zu G\"{o}ttingen, G\"{o}ttingen, 1876,
123-162.




\bibitem{HPV:2010}
P. Hästo, S. Ponnusamy, and M. Vuorinen,
\emph{Starlikeness of the Gaussian hypergeometric
functions}, Complex Var. Elliptic Equ.
\textbf{55} (2010), 173–184.


\bibitem{Kustner:2002}
R.~K\"ustner, \emph{Mapping properties of hypergeometric functions
and convolutions of starlike or convex functions of order $\alpha$},
Comput. Methods Funct. Theory, \textbf{2} (2002), 597--610.

\bibitem{Kustner:2007}
\bysame, \emph{On the order of starlikeness of the shifted
Gauss  hypergeometric function}, J. Math. Anal. Appl. \textbf{334}
(2007), 1363--1385.

\bibitem{LiuPego:2016}
Jian-Guo Liu, R. L. Pego,
\emph{On generating functions of Hausdorff moment sequences},
Trans. Amer. Math. Soc. \textbf{368} (2016), 8499-8518.


\bibitem{MS61}
E.~P. Merkes and W.~T. Scott, \emph{Starlike hypergeometric functions}, Proc.
  Amer. Math. Soc. \textbf{12} (1961), 885--888.


\bibitem{OLBC:2010}
F. W. Olver, D. W. Lozier, R. F. Boisvert, C. W. Clark, \emph{NIMS Handbook of Mathematical
Functions}, Cambridge University Press, 2010.


\bibitem{PonVuor:2001}
S. Ponnusamy, M. Vuorinen,
\emph{Univalence and Convexity Properties for Gaussian Hypergeometric Functions},
Rocky Mountain J. Math. \textbf{31} (2001), 327-353.

 \bibitem{Rus:conv}
S. Ruscheweyh, \emph{Convolutions in {G}eometric {F}unction {T}heory},
  S\'eminaire de Math\'ematiques Sup\'erieures, vol.~83, Les Presses de
  l'Universit\'e de Montr\'eal, Montr\'eal, 1982.


\bibitem{RusSalSu}
S. Ruscheweyh, L. Salinas, and T. Sugawa, \textit{Completely monotone sequences and universally prestarlike functions}, Isr. J. Math. \textbf{171}(2009), 285-304.


\bibitem{RuscheweyhSingh:1986}
S. Ruscheweyh and V. Singh, \textit{On the order of starlikeness of
hypergeometric functions} J. Math. Anal. Appl. \textbf{113} (1986),
1-11.

\bibitem{Silv:1993}
H. Silverman, \emph{Starlike and convexity properties for hypergeometric functions},
J. Math. Anal. Appl. \textbf{172} (1993), 574–581.

\bibitem{SugawaWang:2016}
T.~Sugawa and L.-M.~ Wang, \emph{Note on convex functions of order alpha}, Comput. Methods Funct. Theory, \textbf{16}(2016), 79-92.
%


\bibitem{Wang:2020}
L.-M.Wang, \emph{On the order of convexity for the shifted hypergeometric functions},
 accepted by Comput. Methods Funct. Theory (arXiv:2007.15337).


\bibitem{Umezawa:1952}
T. Umezawa, \emph{Analytic functions convex in one direction},
J. Math. Soc. Japan \textbf{4} (1952),
194-202.

%



\bibitem{Wall:anal}
H. S. Wall, \emph{Analytic Theory of Continued Fractions},
D. Van Nostrand Co. Inc., New York, 1948.
\end{thebibliography}
\end{document}